\documentclass[12pt]{article}

\usepackage{amsmath}
\usepackage{amsthm}
\usepackage{amsfonts}

\parskip\smallskipamount
\newtheorem{theorem}{Theorem}
\newtheorem{lem}{Lemma}
\newtheorem{cor}{Corollary}
\theoremstyle{remark}
\renewcommand{\phi}{\varphi}
\renewcommand{\epsilon}{\varepsilon}

\parskip \smallskipamount

\title{Ergodicity of a stress release point process seismic model with aftershocks}
\author{Pierre Br\'emaud$^{1}$ and Serguei Foss$^{2}$}
\date{}

\begin{document}

\maketitle
\stepcounter{footnote}\footnotetext{INRIA-ENS, D\'epartement d'Informatique, \'Ecole Mormale Sup\'erieure, 45 rue d'Ulm, 75005-Paris, France. Email: pierre.bremaud@ens.fr}
\stepcounter{footnote}\footnotetext{School of Mathematics and Computer Sciences and the Maxwell Institute for Mathematical Sciences, 
  Heriot-Watt University, Edinburgh, Scotland, and 
  Institute of Mathematics, Novosibirsk, Russia. Email: 
S.Foss@hw.ac.uk
The research of S.~Foss was partially supported by the EPSRC Grant EP/E033717/1 and by Framework 7
EURO-NF Grant}
  
\begin{abstract}
We prove ergodicity of a point process earthquake model combining the classical stress release model for primary shocks with the 
Hawkes model for aftershocks. 

{\it Earthquakes, secondary earthquakes, point process, stochastic intensity,
ergodicity, Harris chains, Foster criterion}

\end{abstract}

\section{Introduction}\label{intro}

The times of occurence of earthquakes in a given area of seismic activity form 
a simple point process $N$ on the real line, where $N((a, b])$ is the number of shocks in the time interval $(a, b]$. In the present model, 
the dynamics governing the 
process will be expressed by the stochastic intensity $\lambda (t)$. In intuitive terms (to be precised in the next subsection)
$$ 
\lambda (t) = \lim_{h\downarrow 0} \frac{1}{h} P(N((t, t+h])=1|{\cal F}_t)
$$  
where ${\cal F}_t$ is the sigma-field summarizing the available information at time $t$ (increasing with $t$). 
In the stress release model, for $t\geq 0$, 
$$ 
\lambda (t) = e^{X_0 + ct -\sum_{n=1}^{N((0,t])}Z_n}
$$
where $c>0$ and $\{Z_n\}_{n\geq 1}$ is an i.i.d. sequence of non-negative random variables with finite expectation, 
whereas $X_0$ is some real random variable. 
The process 
$$ 
X(t)=X_0 + ct -\sum_{n=1}^{N((0,t])}Z_n
$$
is ergodic, and the reader is refered to \cite{Last} for a proof and the relevant results concerning 
a generalization of this particular model. 

Another model of interest in seismology is the Hawkes branching process, where the stochastic intensity is 
$$ 
\lambda (t)= \nu (t) + \int_{(0, t]} h(t-s) N(ds), 
$$ 
where $h$ is a non-negative function, called the fertility rate and $\nu$ is a non-negative integrable function. 
Such point process appears in the specialized literature under the name ETAS (Epidemic Type After-Shock; \cite{Ogata}) 
and is used to model the aftershocks (see \cite{DVJ}, 
p. 203). It is well known that the 
corresponding process ``dies out'' in finite time under the condition $\int_0^\infty h(t) \, dt <1$ (\cite{BreNapTor}). 

A model mixing stress release and Hawkes aftershocks is (\cite{Schoenberg})
$$ 
\lambda (t) = e^{X_0 + ct -\sum_{n=1}^{N((0,t])}Z_n} + Y_0 e^{-\alpha t} + k\int_{(0, t]} e^{-\alpha (t-s)} \, N(ds), 
$$ 
where $\alpha >0$. The positive constant $c$ is 
the rate at which the strain builds up. 
If there is a shock at time $t$, then the strain is relieved by the quantity $Z_{N(t)}$. 
Each shock (primary or secondary) at time $t$ generates 
aftershocks according to a Poisson process of intensity $a(s)=k e^{-\alpha (t-s)}$. In this article, we give necessary and sufficient conditions of ergodicity for this model. 
We shall start with a precise mathematical description of it.

\section{Description of the model}\label{Descr}

Let $\varphi : (-\infty ,\infty) \to [0,\infty)$ be a non-decreasing
function such that
$\lim_{x\to -\infty}\varphi (x) = 0$ and $\lim_{x\to\infty} \varphi
(x)
= \infty$. We operate under either one of the following assumptions: 
(a): the function $\varphi$ may be strictly positive everywhere
or (b): it is equal to zero for all $x$ below some level and otherwise strictly increasing.

We are given \\ 
(1) a Poisson field $\Pi$ of intensity 1 in the positive quadrant and \\ 
(2) an i.i.d. family of positive random 
variables $\{ Z_{n}\}_{n\geq 1}$ with a finite
mean,\\
and it is assumed that \\ 
(3) the Poisson field and the i.i.d. family are independent. \\
The above Poisson field and i.i.d. family constitute the probabilistic basis of our model. 

We consider a simple point process $N$ with the following stochastic intensity: 
\begin{equation}\label{spp}
\lambda (t)=\phi 
\left(X_0 +ct -\sum_{n=1}^{N(t)} Z_{n}
\right)+ Y_0e^{-\alpha t} +k \int_{(0,t]}e^{-\alpha (t-s)} \, dN(s), 
\end{equation}
where $N(t):=N((0, t])$, and where $X_0$, $Y_0$, $c$, $k$, and $\alpha$ are as in the introduction. 
This means that the point process is constructed recursively as 
$$ 
N(t)= \int_{(0, t]}\int_R {\mathbf I}(z\leq \lambda (t-)) \Pi (dz \times dz).
$$
Defining 
$$ 
{\cal F}_t:=\sigma \{X_0; Y_0; N(s), 
Z_{N(s)}, 
s \leq t \},
$$
the process $\{\lambda (t)\}_{t\geq 0}$ is then the ${\cal F}_t$-stochastic intensity of $N$ in the sense of \cite{Bre} (see also \cite{DVJ}, \cite{LastBrandt}). \\

In the seismological interpretation, 
\begin{equation}\label{firsts}
\lambda_1 (t)=\phi (X_0 +ct -\sum_{n=1}^{N(t)} Z_{n})
\end{equation}
is the stochastic intensity of the primary shocks, whereas 
\begin{equation}\label{seconds}
\lambda_2 (t)=Y_0e^{-\alpha t} +k \int_{(0,t]}e^{-\alpha (t-s)} \, dN(s) 
\end{equation}
is the stochastic intensity of the aftershocks. \\

The problem is to find a necessary condition for the existence and uniqueness 
of the corresponding stationary  process and, for
any initial distribution of $X_0$ and $Y_0$, of 
the convergence to that distribution, and to prove formally that 
it is also sufficient (under a further smoothnes condition on the distribution
of $Z_i$).

\section{On the ergodicity condition} 

The existence of ergodicity will be
proven in the case 
\begin{equation} 
\frac{k}{\alpha}<1 \label{necessary} 
\end{equation} 
This section shows that this is indeed a natural (intuitive) condition and 
moreover that it is necessary if we seek only those solutions for
which 
the steady-state average intensity 
$\lambda := E[\lambda (t)]$ satisfies $0<\lambda <\infty$. 

We therefore henceforth assume ergodicity.
{F}rom now on we use the notation
$$
X(t) = X_{0}+ct - \sum_{n=1}^{N(t)} Z_{n}
$$
and
$$
Y(t)=Y_0e^{-\alpha t} +k \int_{(0,t]}e^{-\alpha (t-s)} \, dN(s).
$$
The process $(X(t),Y(t)), t\ge 0$ is a time-homogeneous Markov
process with initial value $(X_0,Y_0)$, and
$$
\lambda (t) = \varphi (X(t)) +Y(t).
$$
Further, ergodicity means, in particular, that
there exists a stationary version of the process $(X(t),Y(t))$. 
For such a stationary version, let   
$\lambda_1=E[\varphi (X(t))]$ and $\lambda_{2}= E Y(t)$. Then
$\lambda = \lambda_{1}+\lambda_{2}$, so the finiteness of $\lambda$
implies that of $\lambda_{1}$ and of $\lambda_{2}$. 
Observe that 
\begin{align*} 
E[Y(t)]&= E[Y(0)]e^{-\alpha t} +kE[\int_{(0,t]}e^{-\alpha (t-s)} N(ds)]\\
&= E[Y(0)]e^{-\alpha t} +k\int_0^te^{-\alpha (t-s)} \lambda ds]\\
&= E[Y(0)]e^{-\alpha t} +\lambda \frac{k}{\alpha} (1-e^{-\alpha t})
\end{align*}
where we used Campbell's formula. Therefore, from the
stationarity, 
$E[Y(0)]=\lambda \frac{k}{\alpha}=E[Y(t)]$. Then
$$ 
\lambda = \lambda_1 + \lambda \frac{k}{\alpha} \equiv \lambda_{1}+\lambda_{2}. 
$$   

{\bf The supercritical case.} Suppose, in view of contradiction, that 
$\frac{k}{\alpha}>1$. The last equality then implies that $\lambda =\infty$, which we excluded, or that $\lambda =0$, and then  
$\lambda_1=E[\varphi (X(t))]=0$. Since $\varphi (X(t))\geq 0$, this implies 
$P(\varphi (X(t))= 0)=1$, that is $P(X(t)=-\infty)=1$. Similarly $P(Y(t)=0)=1$. 

{\bf The critical case.} Suppose now, again in view of contradiction, that 
$\frac{k}{\alpha}=1$. The last displayed equality implies then that
$\lambda =\infty$ (excluded) or $\lambda_1=0$ and therefore $P(\varphi
(X(t))= 0)=1$. Then  
$$ 
\lambda (t)=Y(0)e^{-\alpha t} +k\int_{(0,t]}e^{-\alpha (t-s)} N(ds).
$$ 
We show that any point process $N$ with this stochastic intensity and with  finite average intensity is necessarily null (with intensity equal to $0$). Suppose that $\lambda >0$. Clearly, 
\begin{align*} 
P(N(R_+)=0)&=E[P(N(R_+)=0|Y(0))]\\ 
&=E[e^{-\int_0^\infty Y(0)e^{-\alpha t}dt}]\\ 
&\geq e^{-E[Y(0)]\frac{1}{\alpha}} =e^{-\lambda\frac{1}{\alpha}}>0
\end{align*}
and therefore, since we assumed $\lambda <\infty$, we have that $P(N(R_+)=0)>0$. Now, $\{N(R_+)=0\}\subseteq \theta_t\{N(R_+)=0\}=\{N([t,\infty))=0\}$. That is, 
$\{N(R_+)=0\}$ is expanded by the (ergodic) shift, 
and therefore it has probability $0$ or $1$. By the above, this probability must be $1$. We conclude that $\lambda =0$, 
a contradiction. 

Therefore in the critical case there is no solution except the trivial one (no earthquakes).

\section{Explicit expressions for the average rates} 

In this section, we exhibit an interesting feature of the model. We
assume here again ergodicity and the condition $0<\lambda <\infty$. We
continue to consider the model in the stationary regime. 
Writing 
\begin{align*}  
\varphi (X(t))&=\varphi \left(X(0)+t
\left(c -\frac{N(t)}{t}\frac{1}{N(t)}\sum_{n=1}^{N(t)}Z_n
  \right) \right)\\
&=\varphi (X(0)+t(c -\lambda E[Z_1] +\varepsilon (t))), 
\end{align*}
where $\lim_{t\uparrow \infty} \varepsilon (t)=0$ a.s.  
Let $\Delta :=c -\lambda E[Z_1]$. 
Let $\tau$ be the (a.s. finite) random time such that $t\geq \tau$
implies
$|\varepsilon (t)|\leq \frac12 |\Delta|$.

Suppose that $c -\lambda E[Z_1] >0$. We have 
$$ 
E[\varphi (X(t))] \geq E[\varphi (X(t))1_{\{t\geq \tau\}}]\geq E[\varphi (X(0)+t\frac12 \Delta)1_{\{t\geq \tau\}}]
$$
But $\varphi (X(0)+t\frac12 \Delta)1_{\{t\geq \tau\}}\uparrow \infty $
as $t\to\infty$ and therefore  
$$ 
\lambda_1 =E[\varphi (X(t))]\rightarrow \infty 
$$
implying that $\lambda =\infty$ which is excluded. 

Suppose that $c -\lambda E[Z_1] <0$. We show that 
this is imposible. Here $\lim_{t\uparrow \infty} \varphi (X(t))= 0 $ by a similar argument. We prove that  
$\lim_{t\uparrow \infty} E[\varphi (X(t))]= 0 $, $\lambda_1 =0$ 
and therefore $\lambda =0$ which is impossible. 

For the proof that $\lim_{t\uparrow \infty} E[\varphi (X(t))]= 0 $ we 
can make use of the following lemma (in fact taking care of both situations when $c -\lambda E[Z_1] \neq 0$).  



\begin{lem}
If the stationary stochastic process $\{Z(t)\}_{t\geq 0}$ is such that it tends 
almost surely to a deterministic constant $c$ as $t\uparrow\infty$, then it is almost surely equal to this constant. 
\end{lem}

{\sc Proof.} 
Fix $\varepsilon >0$., and consider the set 
$$ 
C=\{\omega ; Z(t,\omega ) \in [c-\varepsilon, c+\varepsilon ] \mbox{ for all } t\geq 0 \}.
$$
Then for all $a>0$, 
$$ 
\theta_a C=\{\omega ; Z(t,\omega ) \in [c-\varepsilon, c+\varepsilon ] \mbox{ for all } t\geq a \}. 
$$
But $\theta_a C \uparrow \Omega$, and therefore $P(C)=P(\theta_a C)\uparrow 1$. So that $P(C)=1$. Since this is true for all $\varepsilon >0$, 
$$ 
P\{Z(t) =c  \}=1, \quad \mbox{for all} \quad t\ge 0.
$$

Therefore, necessarily 
$$
\lambda =\frac{c}{E[Z_1]}. \qquad 
$$

Therefore, in this model, the rate of occurences of earthquakes is given by the physics of stress build up (the constant $c$) and 
stress release ($E[Z_1]$), whereas the 
global rate is shared among primary and secondary earthquakes according to the physics of the aftershocks ($\alpha$ and $k$).

\section{Two embeddings}
We now turn to the technical core of the paper, 
namely the proof of existence of a unique ergodic solution of the model, 
under the condition $\frac{k}{\alpha}<1$ and a further condition on the distribution of
$Z_i$ (see Condition (CZ) in Section 7). The technique used is that of Harris chains, 
and we start as usual by 
studying a natural embedded process. More precisely, 
let $\{(t_n)\}_{n\geq 0}$, with $t_0=0$, be the sequence of time events of $N$, and let for 
each $n\geq 0$, $T_{n+1}:=t_{n+1}-t_n$, $X_n:=X(t_n)$, $Y_n:=Y(t_n)$. We then have 
the recurrence 
equations that exactly reflect the dynamics described in the previous section:
$$ 
X_{n+1}= X_n +cT_{n+1} -Z_{n+1}
$$
and 
$$ 
Y_{n+1}= Y_n e^{-\alpha T_{n+1}} +k
$$
where $S_{n+1}$ is a positive random variable whose hazard rate is, conditionally to $X_0,\ldots, X_n$, $Y_0,\ldots,Y_n$, 
$T_1,\ldots,T_n$ and $Z_1,\ldots, Z_n$
$$ 
\varphi (X_n +cs) + Y_n e^{-\alpha s}.
$$
It is clear that $\{(X_n,Y_n)\}_{n\geq 0}$ is a homogeneous Markov chain. Its transition mechanism is 
fully described by the first transition, which can be implemented as follows:

Let $X_0=x$ and $Y_0=y\ge 0$. On the positive half-plane with the time
$t$ running on the horizontal coordinate axis, draw two
curves: \\
(a) a curve with graph $(t,\varphi (x+ct))$ (that starts from $(0,\varphi (x))$);\\
(b) a curve with graph $(t, -ye^{-\alpha t})$ (that starts from $(0, -y)$).

Consider the projection on the time axis of the Poisson field between the above two curves and let $T_1$ be the 
point of this projection with the smallest $t$-coordinate. It has, as the notation anticipated, 
the required hazard rate $\varphi (x +cs) + y e^{-\alpha s}$. In particular, 
$$
{\mathbf P} (T_1>t \ | \ X_0=x, Y_0=y)
\equiv {\mathbf P}_{x,y} (T>t) =
e^{-\frac{y}{\alpha}
\left(
1-e^{-\alpha t}
\right)}
\cdot
e^{-\int_0^t \varphi (x+cv) dv}.
$$
and 
\begin{eqnarray*}
 X_1&=& x+cT_1 - Z_1,\\
Y_1 &=& ye^{-\alpha T_1} + k.
\end{eqnarray*}

Two lemmas concerning this particular realization of the transition kernel will be useful. 

Let $T_{x,y}$ be a ``generic'' random variable with distribution
$$
{\mathbf P} (T_{x,y} \in \cdot ) = {\mathbf P} (T_1 \in \cdot \ | \ X_0=x, Y_0=y).
$$
Following the comments from above, one can represent $T_{x,y}$ as
\begin{equation}\label{xy}
T_{x,y}= \min \left( T^{(1,x)}, T^{(2,y)} \right)
\end{equation}
where
$T^{(1,x)}$ and $T^{(2,y)}$ are independent and
\begin{eqnarray*}
{\mathbf P} (T^{(1,x)}>t) &=& 
e^{-\int_0^t \varphi (x+cv) dv},\\
{\mathbf P}(T^{(2,y)}>t) &=&
e^{-\frac{y}{\alpha} \left(1-e^{-\alpha t}\right)}.
\end{eqnarray*}

Clearly,
$$
{\mathbf P}(T^{(2,y))}=\infty ) = e^{-y/\alpha} > 0,
$$
for any $y\ge 0$.

\begin{lem}\label{L1}
(1) {F}or any $0\le y_1\le y_2$,
\begin{equation}\label{mon_y}
T^{(2,y_1)}\ge_{st} T^{(2,y_2)}.
\end{equation}
(2)
{F}or any $x_1 < x_2$,
\begin{equation}\label{gege}
T^{(1,x_1)} \ge_{st} T^{(1,x_2)}
\end{equation}
while
\begin{equation}\label{lele}
x_1+cT^{(1,x_1)} \le_{st} x_2+cT^{(1,x_2)}.
\end{equation}
Also, for any $x$, ${\mathbf P} (T^{(1,x)}<\infty ) =1$
and, moreover, for any $a>0$,
$$
{\mathbf P} (T^{(1,x)}>t) e^{at} \to 0,\quad \mbox{as} \quad t\to\infty.
$$
\end{lem}
{\sc Proof.} 
Inequality \eqref{mon_y} is straightforward. 
Inequality \eqref{gege} follows from the monotonicity of
$\varphi$ while inequality \eqref{lele} follows from the 
following coupling construction:\\
Let $t_0>0$ be such that $x_1+ct_0=x_2$. If there is a point of the Poisson
field in $ \{(t,u) \ : 0\le t \le t_0, 0\le u \le \varphi (x_1+ct) \}$,
then  $T^{(1,x_1)} < t_0$ and $x_1+cT^{(1,x_1)}\le x_2$. If however there
is no such a point, then
$$
x_1+cT^{(1,x_1)}= x_2 + cT^{(1,x_2)} \cdot \theta^{t_1}
$$
where $\{\theta^t\}_{t\ge 0}$ is a family of measure-preserving shift transformations.
So, with probability 1
$$
x_1+cT^{(1,x_1)}\le x_2 + cT^{(1,x_2)} \cdot \theta^t =_{st}
x_2+cT^{(1,x_2)}.
$$
The remaining results follow from inequality \eqref{lele} and from the fact that
$\varphi (x) \to\infty$ as $x\to\infty$.

\begin{cor}
For any $x_0$,
\begin{equation}\label{limy}
\sup_{x\ge x_0} {\mathbf E} T_{x,y} \to 0 \quad \mbox{as} \quad y\to\infty.
\end{equation}
Also
\begin{equation}\label{limx}
\sup_{y\ge 0} {\mathbf E} T_{x,y} \to 0 \quad \mbox{as} \quad x\to\infty.
\end{equation}
\end{cor}
{\sc Proof.}
By Lemma 2, $0\le T_{x,y} \le_{st} T_{x_0,y}$ for any $x\ge x_0$.
Clearly, $T^{(2,y)}\to 0$ in probability as $y\to\infty$. Since ${\mathbf E} T^{(1,x_0)}$
is finite, the family of random variables $\{T_{x_0,y}\}_{y\ge 0}$ is uniformly 
integrable, and therefore
$$
\sup_{x\ge x_0} {\mathbf E} T_{x,y} = {\mathbf E} T_{x_0,y} \to 0 \quad
\mbox{as} \quad y\to\infty.
$$
Further, from inequality \eqref{xy}, $T_{x,y}\le_{st} T^{(1,x)}$ where $T^{(1,x)}\to 0$ in
probability. By \eqref{gege} and since ${\mathbf E} T^{(1,0)}$ is finite, 
the family $\{T^{(1,x)}\}_{x\ge 0}$ is uniformly integrable, and therefore ${\mathbf E} T^{(1,x)} \to 0$
as $x\to\infty$, and then \eqref{limx} follows.

\begin{lem}\label{L2}
As $y\to\infty$, 
$$
y{\mathbf E}\left( e^{-\alpha T^{(2,y)}}-1\right) \to -\alpha.
$$ 
\end{lem}
{\sc Proof.}
Indeed,
\begin{eqnarray*}
y{\mathbf E}\left(1- e^{-\alpha T^{(2,y)}}\right)
&=&
y \int_0^1 {\mathbf P} \left(e^{-\alpha T^{(2,y)}}< v \right) dv\\
&=&
y\int_0^1 {\mathbf P} (T^{(2,y)}> \ln v/(-\alpha ) ) dv \\
&=&
y \int_0^1 e^{-\frac{y}{\alpha}(1-v)} dv \quad \left( \mbox{change of variables:} 
\quad u=1-v \right) \\
&=&
y\int_0^1 e^{-\frac{yu}{\alpha}} du \quad \left( \mbox{change of variables:} \
\quad r=\frac{yu}{\alpha} \right)\\
&=& \alpha \int_0^{y/\alpha} e^{-r} dr \rightarrow \alpha,
\end{eqnarray*}
as $y\to\infty$.

{\bf Remark.} As follows from \eqref{xy} and Lemma 2,
if $x\ge 0$, then
$$
T_{x,y} \le_{st} T^{(1,0)}
$$
and, therefore,
\begin{equation}\label{univ1}
\sup_{x>0, y\ge 0}{\mathbf E} T_{x,y} \le {\mathbf E} T^{(1,0)} < \infty.
\end{equation}
One 
may also deduce from Lemma 2 
that 
if $x\le 0$, then
$$
cT_{x,y} \le_{st} |x| + cT^{(1,0)} + Z
$$
where the random variables in the right-hand side are integrable. 
So one can find a universal constant $C>0$ such that
\begin{equation}\label{univ2}
{\mathbf E} T_{x,y} \le C (|x|+1), \quad \mbox{for all}
\quad x\le 0 \quad \mbox{and} \quad y\ge 0.
\end{equation}

Then it follows from \eqref{univ1} and \eqref{univ2} that,
for any negative $x_1$, 
\begin{equation}\label{univ}
\sup_{x> x_1, y\ge 0}{\mathbf E} T_{x,y}  < \infty.
\end{equation}
However the supremum in \eqref{univ} becomes infinite if one replaces
$x_1$ by $-\infty$. 

To keep the supremum finite, we consider a slightly different 
embedding. Again we describe the first transition only.
We fix a sufficiently large positive $v_0$ and a sufficiently
large negative $x_1$ (to be chosen in the next section) and define the
new embedding $\{\widetilde{T}_{x,y}\}$ as follows:\\
(a) if $x\le x_1$, then $\widetilde{T}_{x,y}= \min (T_{x,y},v_0)$ while\\
(b) if $x>x_1$, then $\widetilde{T}_{x,y} = T_{x,y}$.

Then clearly 
\begin{equation}\label{univ_new}
\sup_{x\in (-\infty,\infty ), y\ge 0}{\mathbf E} \widetilde{T}_{x,y}  < \infty.
\end{equation}

Denote by $(\widetilde{X}_n,\widetilde{Y}_n)$
a new time-homogeneous Markov chain obtained by the new embedding. It satisfies
the relations: given $\widetilde{X}_0=x, \widetilde{Y}_0=y$,
if $x>x_1$, then
$$
\widetilde{X}_1=_{st} x+cT_{x,y} - Z 
$$
where $T_{x,y}$ and $Z$ are mutually independent, 
and
$$
\widetilde{Y}_1 =_{st} ye^{-\alpha T_{x,y}} + k,
$$
and if $x\le x_1$, then 
$$
\widetilde{X}_1=_{st} x+c\widetilde{T}_{x,y} - Z {\mathbf I} (T_{x,y}\le v_0)
$$
where $\widetilde{T}_{x,y}$ and $Z$ are mutually independent, and 
$$
\widetilde{Y}_1 =_{st} ye^{-\alpha \widetilde{T}_{x,y}} + k {\mathbf I}
(T_{x,y}\le v_0).
$$

\section{Positive recurrence of the embedded process}

In this section, we show positive recurrence 
of the Markov chain obtained by the second embedding $(\widetilde{X}_n,\widetilde{Y}_n)$
-- see the end of the previous section. We recall some known facts. 

{\bf Definition.} Consider a discrete-time and time-homogeneous Markov
chain $W_{n}, n\ge 0$ on a measurable state space $({\cal W},
{\cal B}_{\cal W})$. A measurable set $V\subseteq {\cal Z}$ is
{\it positive recurrent} if the following two conditions hold:\\
(a) A random variable
$$
\tau_{w}(V) := \min \{ n\ge 1 \ : \ W_{n}\in V \ | \ W_{0}=w\}
$$
is a.s. finite, for all $w\in {\cal W}$;\\
(b)
Moreover,
$$
\sup_{w\in V} E \tau_{w}(V) < \infty.$$

For a Markov chain $ W_{n}$, the following
result is known as {\it Foster's criterion}:

{\bf Proposition 1.}
  Let $L: {\cal W} \to [0,\infty )$
be a measurable function, and let $\widehat{L}$ be a non-negative number.
The set $V=\{ w\in {\cal W} \ : \ L(w) \le \widehat{L}\}$ is positive recurrent
if \\
(i) $\sup_{w\in V} E(L(W_{1}) \ | W_{0}=w )$ is finite;\\
(ii) there exists $\varepsilon > 0$ such that
$$
 E(L(W_{1}) \ | W_{0}=w ) - L(w) \le - \varepsilon,
 $$
for any $w \in {\cal W} \setminus V$.

\begin{theorem}\label{stab}
  Let $k<\alpha $. For suitably chosen $v_0$ and $x_1$, 
there exists a compact set $V$ in $(-\infty,\infty)\times [0,\infty)$ which is positive
recurrent for the Markov chain $(\widetilde{X}_n,\widetilde{Y}_n)$. 
\end{theorem}

{\bf Remark.} A sequence $(X_n,Y_n)$ is a subsequence of $(\widetilde{X}_n,
\widetilde{Y}_n)$.  With arguments similar to those of Theorem \ref{stab}, one can also prove 
 that the same set $V$ is positive
recurrent for the Markov chain $(X_n,Y_n)$.

{\sc Proof of Theorem \ref{stab}.}
We use Foster's criterion, with the following choice of the test function:
$$
 L(x,y) \equiv L_1(x)+L_2(y)= 
 \begin{cases} r_1x+r_2y, & \mbox{ if } x\ge 0 \\
r_3|x|+r_2y, & \mbox{ if } x<0,
\end{cases}
$$
where $r_1,r_2,r_3$ are strictly positive (to be chosen later).

First of all, for any $C_1,C_2>0$,
\begin{equation}\label{CC}
\sup_{|x|\le C_1} \sup_{y\le C_2}
{\mathbf E} (L(\widetilde{X}_1,\widetilde{Y}_1) \ | \ 
\widetilde{X}_0=x, \widetilde{Y}_0=y) < \infty.
\end{equation}
Indeed, let $r=\max (r_1,r_2,r_3)$. Then, for any $(x,y)$ from the set
above,
$$
\widetilde{T}_{x,y} \le_{st} T^{(1,-C_1)} :=\widehat{T},
$$
and $\widehat{T}$ has a finite mean. 
Therefore, for all $(x,y)$ from this set,
$$
{\mathbf E} (L(\widetilde{X}_1,\widetilde{Y}_1) \ | \ 
\widetilde{X}_0=x, \widetilde{Y}_0=y) \le
r(C_1+C_2+{\mathbf E} \widehat{T} + {\mathbf E} Z + k) < \infty,
$$
and \eqref{CC} follows. 

Now we impose several constraints on the coefficients $r_1,r_2,r_3$. 
First, we assume that
\begin{equation}\label{first1}
 r_3 < r_1 
\end{equation}
and 
\begin{equation}\label{second1}
r_1{\mathbf E} Z > r_2k.
\end{equation}
Let $\alpha -k=2\Delta >0$. 
We also assume that 
\begin{equation}\label{third1}
r_2\Delta > r_3 {\mathbf E}Z.
\end{equation}
In the proof, we use only conditions \eqref{first1} - \eqref{third1}
which are, in particular, satisfied if $r_1 >> r_2 >> r_3 >0$.

Now we proceed to show that all the differences
$$
E (L(\widetilde{X}_{1},\widetilde{Y}_{1}) \ | \ 
(\widetilde{X}_{0},\widetilde{Y}_{0})=(x,y)) - L(x,y)
$$
are bounded above by some negative constant if 
$(x,y)$ is outside the set $[x_1,x_0] \times [0,y_0]$ where
$x_1,x_0$ and $y_0$ will be chosen in the proof.  
For this, consider separately two cases: (a) $x>0$ and (b) $x\le 0$.

\vspace{0.3cm}

{\large \bf Case $x>0$}

In this case, the one-step embedding is the natural one, so
we may write $(X_1,Y_1)$ instead of $(\widetilde{X}_1,
\widetilde{Y}_1)$.

Let 
\begin{equation}\label{eq1}
3\gamma = \min \{ r_2\Delta - r_3 {\mathbf E}Z,
r_1 {\mathbf E} Z - r_2 k \} >0.
\end{equation}

Choose  $x_0 >0$ so big that 
\begin{equation}\label{eq2}
r_1c {\mathbf E} T^{(1,x_0)} \le \gamma,
\end{equation}
and
\begin{equation}\label{eq3}
(r_1+r_3) {\mathbf E} (Z-x_0)^+ \le \gamma.
\end{equation}
By Lemma 3, we may choose $y_0>0$ so large that
\begin{equation}\label{eq4}
y{\mathbf E} \left(1-e^{-\alpha T^{(2,y)}}\right) - k \ge \frac{5}{3}\Delta, \quad
\mbox{for all} \quad y\ge y_0
\end{equation}
and that
\begin{equation}\label{eq5}
r_1c{\mathbf E} T_{0,y_0} \le \gamma.
\end{equation}
Then choose $v_0>0$ so large that
\begin{equation}\label{eq6}
y{\mathbf E} \left(1-e^{-\alpha \min (v_0,T^{(2,y)})}\right) - k \ge \frac{4}{3}\Delta, \quad
\mbox{for all} \quad y\ge y_0
\end{equation}
and that the following inequality holds:
\begin{equation}\label{eq10}
\frac{r_3cv_0}{2}e^{-y_0/\alpha} > \gamma + r_3 {\mathbf E} Z + r_2k.
\end{equation}

Write for short
$$
{\mathbf E}_{x,y} L(X_1,Y_1)= {\mathbf E} (L(X_1,Y_1) \ | \ X_0=x, Y_0=y)
= {\mathbf E}_{x,y} L_1(X_1) + {\mathbf E}_{x,y} L_2(Y_1).
$$
If $x>0$, then
\begin{eqnarray*}
{\mathbf E}_{x,y}L_1(X_1) - L_1(x)
&=&
r_1 {\mathbf E} \left((x+cT_{x,y} -Z) {\mathbf I} (x+cT_{x,y} -Z>0)
\right)\\
&+& 
r_3 {\mathbf E} \left((-x-cT_{x,y} +Z) {\mathbf I} (x+cT_{x,y} -Z\le 0)
\right) - r_1x\\
&= &
c{\mathbf E}\left(T_{x,y} \left(
r_1{\bf I} (x+cT_{x,y}-Z>0) - r_3 {\bf I} (x+cT_{x,y}-Z\le 0) \right)\right) \\
&+&
r_1 {\mathbf E} (x-Z) + (r_1+r_3) {\mathbf E} \left(
(-x+Z) {\mathbf I} (x+cT_{x,y}-Z\le 0) \right) - r_1x \\
&\le &
r_1 c {\mathbf E} T_{x,y} - r_1{\mathbf E} Z + (r_1+r_3)
{\mathbf E} (Z-x)^+
\end{eqnarray*}
and
$$
{\mathbf E} L_2(Y_1) - L_2(y) =
r_2 (y{\mathbf E}e^{-\alpha T_{x,y}}+k) -r_2y.
$$
In particular, if $x\ge x_0$, 
$$
{\mathbf E}_{x,y}L_1(X_1) - L_1(x)
\le 
r_1 c {\mathbf E} T^{(1,x_0)} - r_1 {\mathbf E} Z + (r_1+r_3) {\mathbf E} (Z-x_0)^+ 
$$
(where we used representation \eqref{xy} and  Lemma 2) and
\begin{eqnarray*}
{\mathbf E} L_2(Y_1) - L_2(y) \le r_2k,
\end{eqnarray*}
so in view of \eqref{eq1}, \eqref{eq2}, and \eqref{eq3},
$$ 
{\mathbf E}_{x,y} L(X_1,Y_1) - L(x,y)\le \gamma + \gamma - 3\gamma = -\gamma.
$$
Furthermore, if $y\ge y_0$ and $0\le x \le x_0$, then, by Lemma 2,
$$
{\mathbf E}_{x,y} L_1(X_1) - L_1(x) \le  
r_1c {\mathbf E}T_{0,y_0} + r_3 {\mathbf E}Z
$$
and, by inequality \eqref{eq4},
$$
{\mathbf E} L_2(Y_1) - L_2(y) \le -r_2 \Delta,
$$
so
$$
{\mathbf E}_{x,y} L(X_1,X_2) - L(x,y)
\le r_1c {\mathbf E}T_{0,y_0} + r_3 {\mathbf E}Z - r_2 \Delta \le 
\gamma - 3\gamma \le -\gamma,
$$
by \eqref{eq1} and \eqref{eq5}. 
So if $x>0$ and if either $x\ge x_0$ or $y\ge y_0$, 
then
$$
E_{x,y} L(X_{1},Y_{1}) - L(x,y) \le - \gamma.
    $$


\vspace{0.3cm}
{\large \bf Case $x\le 0$}

{F}or the time being, fix any value of $x_1 <0$. 
{F}irst, we observe that if $y\ge y_0$ where $y_0$ satisfies inequalities
\eqref{eq4} and \eqref{eq5}, then again the increments
$E_{x,y} L(X_{1},Y_{1}) - L(x,y)$ have 
a ``uniformly''
negative drift
in all $ x_1 < x \le 0$. 
Indeed, if $\widetilde{X}_0=x \in (x_1, 0]$, then, for any $y\ge 0$, 
$\widetilde{X}_1 =x+ T_{x,y}-Z_1$ admits the following bounds:
\begin{equation}\label{lower0}
|\widetilde{X}_1| {\mathbf I} (\widetilde{X}_1\le 0) \le_{a.s.}
(|x|+ Z_1) {\mathbf I} (\widetilde{X}_1\le 0)
\end{equation}
and
\begin{equation}\label{upper0}
\widetilde{X}_1 {\mathbf I} (\widetilde{X}_1>0) = \max (0, 
\widetilde{X}_1)
= 
\max (0, \min (x+cT^{(1,x)}, x+T^{(2,y)})).
\end{equation}
From Lemma 2 and from independence of $T^{(1,x)}$ and $T^{(2,y)}$, we obtain
\begin{equation}\label{upper}
\widetilde{X}_1 {\mathbf I} (
\widetilde{X}_1>0)\le_{st} \min (c {T}^{(1,0)}, cT^{(2,y)})
= c{T}_{0,y}.
\end{equation}
Therefore, for any $x_1 <x\le 0$
\begin{eqnarray*}
{\mathbf E}_{x,y} |L_1(\widetilde{X}_1)| 
&=&
r_1 {\mathbf E}_{x,y} (\widetilde{X}_1 {\mathbf I} (\widetilde{X}_1> 0))+
r_3 {\mathbf E}_{x,y} (|\widetilde{X}_1| {\mathbf I}(
\widetilde{X}_1\le 0)) \\
 &\le &
r_1 
{\mathbf E}_{x,y} (|\widetilde{X}_1|\cdot {\bf I} (
\widetilde{X}_1>0))
+ r_3 {\mathbf E}_{x,y} ((|x|+{\mathbf E} Z_1){\mathbf I} 
(\widetilde{X}_{1}\le 0))\\
&\le &
r_1c{\mathbf E} {T}_{0,y} + r_3|x| + r_3{\mathbf E} Z.
\end{eqnarray*}
Since $y_0$ satisfies unequalities \eqref{eq4} and \eqref{eq5},
we have, for all $x_1< x\le 0$ and $y\ge y_0$,
$$
{\mathbf E}_{x,y} L(\widetilde{X}_1,\widetilde{Y}_1)-L(x,y) \le 
r_1c{\mathbf E} {T}_{0,y_0} + r_3{\mathbf E} Z - r_2 \Delta
\le - \gamma.
$$

We now choose $x_1 <<-1$ so large that the increment
of $L(\widetilde{X},\widetilde{Y})$ has a uniformly negative drift
for all $x\le x_1$.
We start with the assumption that
\begin {equation}\label{x1x1}
x_1 \le -cv_0.
\end{equation}
Therefore, if $\widetilde{X}_0=x\le x_1$ and $\widetilde{Y}_0=y$, then 
$$
L_1(\widetilde{X}_1) = r_3 (-x-c\widetilde{T}_{x,y} +
Z {\mathbf I} (T_{x,y}\le v_0)) \le r_3 (-x-c\widetilde{T}_{x,y} +Z)
$$
and
$$
L_2(\widetilde{Y}_1) =
r_2 \left(
ye^{-\alpha \widetilde{T}_{x,y}}+k {\mathbf I} (T_{x,y} \le v_0)
\right)
\le
r_2
\left(
ye^{-\alpha \widetilde{T}_{x,y}}+k \right). 
$$
We impose two additional constraints on $x_1$ making it 
even more negatively large. 
Since $\varphi (x) \to 0$ as $x\to -\infty$, 
one can choose $x_1 < -cv_0$, $x_1 << -1$ such that 
\begin{equation}\label{choiceofx1}
e^{-\int_0^{v_0} \varphi (x_1+cv) dv }\ge 1/2.
\end{equation}
Secondly, 
it follows that
 $T^{(1,x)}\to \infty$ in probability as $x\to -\infty$, and therefore, from
\eqref{eq6}, one can choose $x_1 \le -cv_0$ such that
\begin{equation}\label{eq7}
y{\mathbf E} \left(1-e^{-\alpha \widetilde{T}_{x,y}}\right) - k \ge \Delta, \quad
\mbox{for all} \quad y\ge y_0 \quad \mbox{and} \quad x\le x_1.
\end{equation}

Assume that $x_1$ satisfies all of the three conditions
\eqref{x1x1}--\eqref{eq7}. 
If $y\ge y_0$ then, for any $x\le x_1$,
\begin{eqnarray*}
{\mathbf E}_{x,y} L(\widetilde{X}_1,\widetilde{Y}_1) - L(x,y)
&\le &
r_3 (-c{\mathbf E}\widetilde{T}_{x,y} +{\mathbf E}Z)
+
r_2
\left(
y{\mathbf E}e^{-\alpha \widetilde{T}_{x,y}}+k \right) - r_2y\\
&\le &
-r_2 \Delta + r_3 {\mathbf E}Z \le -\gamma,
\end{eqnarray*}
by \eqref{eq7}.

If instead $y\le y_0$ and $x\le x_1$, then
$$
{\mathbf P} (T_{x,y}>v_0) \ge 
{\mathbf P} (T_{x,y_0}>v_0) \ge \frac{1}{2} e^{-y_0/\alpha},
$$
since the random variables $T_{x,y}$ are stochastically decreasing in $y$
(again by Lemma 2).
Therefore, for
$
\widetilde{T}_{x,y}= \min (v_0, T_{x,y}),
$
\begin{equation}\label{defK}
{\mathbf E} T_{x,y} \ge 
{\mathbf E}
\widetilde{T}_{x,y} \ge v_0 \cdot {\mathbf P} (T_{x,y}\ge v_0)
\ge \frac{v_0}{2}e^{-y_0/\alpha} ,
\end{equation}
and
\begin{eqnarray*}
{\mathbf E}_{x,y} L(\widetilde{X}_1,\widetilde{Y}_1) - L(x,y)
&\le &
r_3 (-c{\mathbf E}\widetilde{T}_{x,y} +{\mathbf E}Z)
+
r_2
\left(
y{\mathbf E}e^{-\alpha \widetilde{T}_{x,y}}+k \right) - r_2y\\
&\le &
-\frac{r_3cv_0}{2}e^{-y_0/\alpha} + r_3{\mathbf E} Z + r_2k < -\gamma,
\end{eqnarray*}
due to \eqref{eq10}.

As an outcome, we have that if $y_0$ satisfies conditions \eqref{eq4}--\eqref{eq5},
if $v_0$ satisfies \eqref{eq6}--\eqref{eq10}, and if $x_1$
satisfies \eqref{x1x1}--\eqref{eq7},
then the increments of ${\mathbf E} L(\widetilde{X}_1,
\widetilde{Y}_1)$ have a drift bounded above by $-\gamma$ for all
initial values such that either $x\le x_1$, or $x\le 0$ and $y\ge y_0$.  

The set 
$$
V= [x_1,x_0]\times [0,y_0]
$$ 
is therefore positive
recurrent for the Markov chain $(\widetilde{X}_n,\widetilde{Y}_n)$.

Also, as follows from
the classical proof of Foster's criterion,  
for any initial value $(\widetilde{X}_0,
\widetilde{Y}_0)=(x,y)$, a random variable
$$
\tau_{x,y} (V) = \min \{ n\ge 1 \ : \ 
(\widetilde{X}_n,\widetilde{Y}_n)\in V \ | \
(X_0,Y_0)=(x,y) \}
$$
is almost-surely finite and, moreover, 
there exists an absolute constant
$C>0$ such that
$$
{\mathbf E} \tau_{x,y}(V) \le C(L(x,y)+1),
$$
for all $(x,y)$ (see, e.g., \cite{MT} or \cite{FK}).
The proof of Theorem \ref{stab} is complete.

\section{Harris ergodicity}


We recall the following classical result (see for instance \cite{MT}).

{\bf Proposition 2.}
Assume that a Markov chain $W_n, n\ge 0$ taking valued in a measurable space
$({\cal W, B_W})$ is aperiodic and that there exists a positive
recurrent set ${V}$ that admits a minorant measure, i.e. there exist
%
a positive integer $m$, a positive
$p\le 1$ and a probability measure $\mu$ 
such that
\begin{equation}\label{minor}
{\mathbf P} (W_m\in \cdot \ | W_0=w\in {{V}}) \ge p \mu (\cdot ).
\end{equation}
Then the Markov chain is {\it Harris ergodic}, which means that there exists a unique stationary
distribution (say $\pi$) and that, for any initial value $W_0=w$, there is a convergence
of the distributions of $W_n$ to the stationary one in the {\it total variation norm},
$$
\sup_{B\in {\cal B}_{\cal W}} |P(W_n\in B) - \pi (B)| \to 0, \quad n\to \infty.
$$

In practice, the most technical part in applying this criterion is to verify the 
aperiodicity.
There are a number of sufficient conditions available for the Markov chain to be aperiodic and
Harris ergodic.

We mention two of them. The most common is the following condition.

{\bf Sufficient condition 1 (SC1).} 
A Markov chain $W_n$ is Harris ergodic if there exists 
a positive recurrent set ${V}$ such that
condition \eqref{minor} holds with $m=1$ and with $\mu$ such that $\mu ({V})>0$.

However, in our proof, it seems to be easier to verify another --- slightly more general --- sufficient condition.

{\bf Sufficient condition 2 (SC2)}
A Markov chain $W_n$ is Harris ergodic if there exists 
a positive recurrent set ${V}$ such that
condition \eqref{minor} holds with a finite number of different values of $m$, say
$m_i$, $i=1,2,\ldots , k$ which are such that
$$
g.c.d. \{ m_i, 1\le i \le k \} =1.
$$

We will apply condition (SC2) with $k=2$ and with $m_1=2$ and $m_2=3$.
{F}or that, we introduce a condition on the distribution of $Z$ which leads to (SC2).

{\bf Condition (CZ).} There exist $0\le z_1 < z_2<\infty$ such that, for some $h>0$
and for any $[u_1,u_2]\subseteq [z_1,z_2]$,
$$
{\mathbf P} (Z\in [u_1,u_2]) \ge h(u_2-u_1).
$$
In other words, the distribution of $Z$ has an absolutely continuous (w.r. to Lebesgue
measure) component whose density function is above level $h$ everywhere in the
interval $[z_1,z_2]$.

\begin{theorem}\label{th2}
Assume condition (CZ) to hold. Then the Markov chain $(X_n,Y_n)$ is Harris ergodic.
\end{theorem}

{\sc Proof.}
We may assume without loss of generality that $z_2-z_1 \le x_0-x_1$. 

Let $\widetilde{y}_1=y_1+k$ and 
$$
\widetilde{x}_0= \inf \{ x\ge x_0 \ :  \ \varphi (x) >0\} + z_2+2,
$$
and let
$$
V_1 = [x_1,\widetilde{x}_0]\times [0,\widetilde{y}_0],
$$
so $V \subset V_1$. 
Then, for any $(x,y)\in V$, ${\mathbf P}_{x,y} (Y_1 \le \widetilde{y}_0)=1$, so,
by Lemmas 2--3, 
\begin{eqnarray*}
{\mathbf P}_{x,y} ((X_1,Y_1)\in V_1)
&=&
{\mathbf P}_{x,y} (X_1 \in [x_1,\widetilde{x}_0])\\
&\ge &
{\mathbf P} (x+T_{x,y} \in [\widetilde{x}_0+z_1-1,\widetilde{x}_0+z_1], Z_1 \in [z_1,z_2] )
\\
&\ge &
{\mathbf P} (x_1+cT^{(1,x_1)}\in [\widetilde{x}_0+z_1-1, \widetilde{x}_0+z_1]) \\
&\times &
{\mathbf P} (cT^{(2,\widetilde{y}_0)}> \widetilde{x}_0+z_1-x_1)
{\mathbf P} (Z_1 \in [z_1.z_2]).
\end{eqnarray*}
Denote by $R_0$ the value of the rightmost side of the above inequality (note that it is positive). Then, for
any $(x,y)\in V_1$, 
$$
{\mathbf P}_{x,y} ((X_1,Y_1) \in V_1) \ge R_0>0.
$$

Take some small positive $\varepsilon <(z_2-z_1)/4$ (to be specified later). 
Choose $t_2>0$ so large that 
$x_2:=x_1+ct_2 > \widetilde{x}_0+z_2$ and 
$y_0e^{-\alpha (x_2-\widetilde{x}_0)}\le \varepsilon$.
Let $b=\varphi (x_2-z_2)$ and note that $b>0$.
Then, for any $(x,y)\in V_1$, 
$$
{\mathbf P} (x +cT_{x,y} \in [x_2,x_2+\varepsilon ], ye^{-\alpha T_{x,y}}\le
\varepsilon ) 
\ge 
b\varepsilon {\mathbf P} (T^{(1,x_1)}>(x_2-x_1)/c)
{\mathbf P}( T^{(2,y_0)}> (x_2-x_1)/c). 
$$
Denote by $R_1$ the right-hand side of the inequality above (which is a positive number). 
Then, for any $(x,y)\in V_1$ and for $(X_1,Y_1)=(x+T_{x,y}-Z_1, ye^{-\alpha T_{x,y}}+k)$,
$$
{\mathbf P}_{x,y} ((X_1,Y_1)\in [x_2-z_1-\varepsilon, x_2-z_1]\times [k,k+\varepsilon ])
\ge R_1 \frac{b\varepsilon}{z_2-z_1} =: R_2 >0.
$$
Let
$$
\widehat{V} = [x_2-z_1-\varepsilon, x_2-z_1]\times [k,k+\varepsilon ].
$$
{F}rom the construction above, one may conclude that, for any $(x,y)\in V$,
\begin{equation}\label{ffirst}
\inf_{(x,y)\in V} {\mathbf P}_{x,y} ((X_1,Y_1)\in \widehat{V})\ge R_2 >0
\end{equation} 
(since $V \subset V_1$) and then that, by the Markov property,
\begin{equation}\label{fsecond}
{\mathbf P}_{x,y} ((X_2,Y_2)\in \widehat{V}) \ge
R_0 \cdot \inf_{(x,y)\in V_1} {\mathbf P}_{x,y} ((X_1,Y_1)\widehat{V}) = R_0R_2 >0.
\end{equation}

Now take $\varepsilon >0$ so small that one can choose positive numbers $t_3$ and $t_4$
such that $t_4>t_3>z_2$, 
that
$$
k_2:= ke^{-\alpha t_3} > (k+\varepsilon ) e^{-\alpha t_4} =: k_1,  
$$
and that
$$
\delta := \varepsilon + c(t_4-t_3) < \frac{z_2-z_1}{2}.
$$
Then, for any $y\in [k,k+\varepsilon ]$, we have the inclusion 
$[k_1,k_2]\subseteq [ye^{-\alpha t_4}, ye^{-\alpha t_3}]$. 

For any $(x,y)\in \widehat{V}$ denote by $g_{x,y}(u)$ a density function of random variable
$ye^{-alpha T_{x,y}}$ (which clearly has an absolutely continuous distribution).


Then direct computations show that  
\begin{equation}\label{ccc}
c_0:=\inf_{(x,y)\in \widehat{V}} \inf_{u\in [k_1,k_2]} g_{x,y}(u).
\end{equation}
Indeed, let 
$$
c_1 = \inf_{0\le t \le k_2/k_1} \frac{\ln (1+t)}{t} \quad \mbox{and} \quad 
c_2 = \inf_{(x,y)\in \widehat{V}} 
\inf_{\ln (k/k_2)\le a < b\le \ln ((k+\varepsilon )/k_1)}
\frac{{\mathbf P} (a\le T_{x,y}\le b)}{b-a}.
$$
Then both $c_1$ and $c_2$ are positive and, for $[a,b]\subseteq [k_1,k_2]$,
$$
{\mathbf P}
(ye^{-\alpha T_{x,y}}\in [a,b]) = {\mathbf P}
\left(
\frac{\ln (y/b)}{\alpha} \le T_{x,y} \le \frac{\ln (y/a)}{\alpha} \right)
\ge \frac{c_1c_2}{\alpha} (b-a),
$$
so \eqref{ccc} holds with $c_0=c_1c_2/\alpha$.

Furthermore, let $x_3 = x_2-z_1-\varepsilon + t_3$. Note that if $(x,y)\in \widehat{V}$
and $T_{x,y}\in [t_3,t_4]$, then 
$x+cT_{x,y} \in [x_3,x_3+\delta]$. Then, by condition (CZ), given 
$x+cT_{x,y}=v \in [x_3,x_3+\delta]$, the random variable $v-Z_1$ has an absolutely
continuous component with a uniform 
distribution on the interval $[x_3-(z_1+z_2)/2, x_3-z_1]$.

We may therefore conclude that, for any $(x,y) \in \widehat{V}$,
\begin{equation}\label{fthird}
{\mathbf P}_{x,y} ((X_1,Y_1) \in \cdot )\ge
2(z_2-z_1)^{-1}hc_0^{-1} \mu (\cdot )
\end{equation}
where $\mu$ is a two-dimentional uniform distribution on the rectangle
$V_2 :=[x_3-(z_1+z_2)/2, x_3-z_1]\times [k_1,k_2]$, coefficient $h$ is from 
condition (CZ), and $c_0$ is from \eqref{ccc}.

It follows from inequalities \eqref{ffirst}, \eqref{fsecond}, and \eqref{fthird}, 
that condition (SC2) is satisfied with $k=2$ and with $m_1=2$ and $m_2=3$, and this completes the proof.

\begin{cor} Assume again that $k<\alpha$ and that condition (SZ) holds. 
Consider the Markov chain $(\widetilde{X}_n,\widetilde{Y}_n)$ and let
$0\le S_1 < S_2 < \ldots S_k < \ldots$ be consequtive times of the ends of 
discrete-time ``cycles''
where random vectors 
$(\widetilde{X}_{S_k}, \widetilde{Y}_k)$ 
have uniform distribution in the rectangle $V_2$.
Then random variables $l_k=S_k-S_{k-1}$, $k\ge 2$ are i.i.d. with a finite mean.
\end{cor}

A proof of this result can be found, for instance, in \cite{Asm}.

\section{Stability in continuous time}

\begin{theorem}\label{stab_cont}
Under condition (CZ), \\
(1) there exists a unique stationary version of the continuous-time Markov process
$(X(t),Y(t)$ (which is also ergodic);\\
(2)  for any initial value $X(0)=X_0=x$, $Y(0)=Y_0=y$, the 
process $(X(t),Y(t))$ converges to the stationary one in the total variation norm. 
\end{theorem}

{\sc Proof.}
Consider again the embedded Markov chain $(\widetilde{X}_n,\widetilde{Y}_n)$ and
its cycles of length $l_k$. Then the corresponding cycles in continuous time are
defined as
$$
L_k = \sum_{i=S_{k-1}+1}^{S_k} \widehat{T}_{\widehat{X}_i,\widehat{Y}_i}, 
\quad k=1,2,\ldots
$$
which are again i.i.d. for $k\ge 2$.

Then a proof of the
theorem follows from the two results below, Statement 1 and Statement 2
(see, for instance, \cite{Asm}, Proposition 3.8, p.203 or \cite{BF}, Section 7 or \cite{Bor}, Chapter 3).

{\bf Statement 1.} The distribution of random variables $L_k$, $k\ge 2$ 
has an absolutely continuous component  with a density function which is separated from
$0$ on a finite time interval of positive length.

{\bf Statement 2.} ${\mathbf E} L_2$ is finite.

Statement 1 may be verified directly using arguments similar to those in 
the previous section.
Furthermore, since 
$
C:= \sup_{x,y} {\mathbf E} \widetilde{T}_{x,y} < \infty , 
$

$$
{\mathbf E} L_2 \le C {\mathbf E} l_2 < \infty ,
$$ 
and the result follows.



\end{document}